\documentclass[11pt,titlepage]{amsart}
\usepackage{amsmath,amsthm, amscd, amssymb, amsfonts}

\usepackage[all]{xy}
\usepackage{mathrsfs}
\usepackage{enumitem}
\usepackage[T2A,T1]{fontenc}

\usepackage[all]{xy}

\usepackage[ansinew]{inputenc}
\usepackage{graphicx,fancyhdr}

\usepackage[dvips, dvipsnames, usenames]{color}

\newcommand{\diag}{\operatorname{diag}}

\newcommand{\comment}[1]{}

\newcommand{\zh}{\mbox{\usefont{T2A}{\rmdefault}{m}{n}\cyrzh}}
\newcommand{\Zh}{\mbox{\usefont{T2A}{\rmdefault}{m}{n}\CYRZH}}

\newcommand{\supp}{\operatorname{supp}}

\numberwithin{equation}{section}
\newtheorem{theorem}{Theorem}[section]
\newtheorem{lemma}[theorem]{Lemma}
\newtheorem{coro}[theorem]{Corollary}
\newtheorem{conjecture}[theorem]{Conjecture}
\newtheorem{prop}[theorem]{Proposition}

\newtheorem{paso}{Step}

\theoremstyle{definition}
\newtheorem{definition}[theorem]{Definition}
\newtheorem{example}[theorem]{Example}
\newtheorem{question}{Question}

\theoremstyle{remark}
\newtheorem{remark}[theorem]{Remark}

\newcommand{\pf}{\begin{proof}}
  \newcommand{\epf}{\end{proof}}

\newcommand{\fjdos}{\hspace{-1pt}j + \frac{1}{2}}
\newcommand{\fkdos}{\hspace{-1pt}k + \frac{1}{2}}

\newcommand{\inc}{\mathscr{I}}

\newcommand{\ba}{ \mathbf{a}}

\newcommand{\ku}{ \Bbbk}
\newcommand{\kut}{ \ku^{\times}}

\newcommand{\G}{\mathbb G}

\newcommand{\ghost}{\mathscr{G}}

\newcommand{\I}{\mathbb I}
\newcommand{\Iw}{\mathbb I^{\dagger}}

\newcommand{\J}{\mathbb J}
\newcommand{\N}{\mathbb N}
\newcommand{\bp}{\mathbf{p}}
\newcommand{\bq}{\mathbf{q}}

\newcommand{\Q}{\mathbb Q}
\newcommand{\Sb}{\mathbb S}
\newcommand{\Uu}{\mathbb U}
\newcommand{\V}{\mathbb V}

\newcommand{\Z}{\mathbb Z}

\newcommand{\D}{\mathcal{D}}

\newcommand{\cI}{\mathcal{I}}
\newcommand{\cJ}{\mathcal{J}}

\newcommand{\cR}{\mathcal{R}}

\newcommand{\T}{\mathcal{T}}
\newcommand{\cV}{\mathcal{V}}
\newcommand{\X}{\mathcal{X}}
\newcommand{\Xx}{\mathscr{X}}
\newcommand{\Xf}{\X_{\text{fin}}}
\newcommand{\Xif}{\X_{\infty}}

\newcommand{\Fg}{\mathfrak F}

\newcommand{\lstr}{\mathfrak L}
\newcommand{\cyc}{\mathfrak C}
\newcommand{\pos}{\mathfrak P}

\newcommand{\End}{\operatorname{End}}
\newcommand{\id}{\operatorname{id}}
\newcommand{\gr}{\operatorname{gr}}
\newcommand{\GK}{\operatorname{GKdim}}

\newcommand{\ydh}{{}^{H}_{H}\mathcal{YD}}

\newcommand{\toba}{\mathscr{B}}
\newcommand{\ot}{\otimes}

\newcommand{\ydG}{{}^{\ku \Gamma }_{\ku \Gamma }\mathcal{YD}}

\newcounter{tabla}\stepcounter{tabla}

\begin{document}
  
  \title[Nichols algebras with finite Gelfand-Kirillov dimension]{On Nichols algebras of infinite rank with finite Gelfand-Kirillov dimension}
  
\author[Andruskiewitsch, Angiono, Heckenberger]
{Nicol\'as Andruskiewitsch, Iv\'an Angiono, Istv\'an Heckenberger}
  
\thanks{\noindent\hspace{-12pt} N. A. (corresponding author): FaMAF-CIEM (CONICET), Universidad Nacional de C\'ordoba,
 Medina A\-llen\-de s/n, Ciudad Universitaria, 5000 C\' ordoba, 
 Rep\'ublica Argentina. \newline\texttt{andrus@famaf.unc.edu.ar}
 \newline \newline
I. A.: FaMAF-CIEM (CONICET), Universidad Nacional de C\'ordoba,
 Medina A\-llen\-de s/n, Ciudad Universitaria, 5000 C\' ordoba, 
 Rep\'ublica Argentina. \newline\texttt{angiono@famaf.unc.edu.ar}
 \newline \newline
I. H.: Philipps-Universität Marburg,
 	Fachbereich Mathematik und Informatik,
 	Hans-Meerwein-Straße,
 	D-35032 Marburg, Germany. \newline\texttt{heckenberger@mathematik.uni-marburg.de}
 The work of N. A. and I. A.  was partially supported by CONICET,
 Secyt (UNC), the MathAmSud project GR2HOPF. The work of I. A. was partially supported by ANPCyT (Foncyt).
 The work of N. A., respectively I. A., was partially done during a visit to the University of Marburg, 
 respectively the MPI (Bonn), supported by the Alexander von Humboldt Foundation. }

  \begin{abstract}
 We classify infinite-dimensional decomposable braided vector spaces arising from abelian groups whose components are either points or blocks
 such that the corresponding Nichols algebras have finite Gelfand-Kirillov dimension. In particular we exhibit examples where  the  Gelfand-Kirillov dimension attains any natural number.
  \end{abstract}
  
  \maketitle


\section{Introduction}\label{section:introduction}

The study of  Hopf algebras with finite Gelfand-Kirillov dimension (abbreviated $\GK$)
received considerable attention in the last years, see e.g. 
\cite{triang,B-seattle,B-turkish,BG,BGLZ,EG,G-survey,R-quantum-groups}
and references therein, or \cite{A-icm} for $\GK =0$. 
By the lifting method \cite{AS Pointed HA}, one is naturally led to consider the problem of classifying Nichols algebras 
with finite Gelfand-Kirillov dimension.

Let $\ku$ be an algebraically closed field of characteristic 0,
let $\Gamma$ be an abelian group and let $\V\in \ydG$ such that  $\dim \V$ is infinite and countable. 
In this paper we contribute to the following question: when $\GK \toba(\V) < \infty$?

We first show that the underlying braided vector space is locally finite, see Theorem \ref{prop:locally-finite-Z}. 
Then we consider two classes of braided vector spaces with infinite basis: those of diagonal type, and those that are sums of points and blocks.
The classification of those $\V$ of diagonal type with connected diagram
 such that $\GK \toba(\V) = 0$ follows a well-known pattern, see Proposition \ref{prop:gk0}. 
We have proposed in \cite[1.2]{triang}:

\begin{conjecture}\label{conjecture:nichols-diagonal-finite-gkd} 
If $V$ is a finite-dimensional braided vector space of diagonal type such that $\GK \toba(V) < \infty$, then it has an arithmetic root system.
\end{conjecture}
  
In other words, such $V$ should belong to the classification in \cite{H-classif}. 
Some evidence on the validity of this Conjecture is offered in \cite{diag}, where we show that it is valid for rank 2 and for affine Cartan type.
Assuming this Conjecture, it is not difficult to prove:

\begin{prop}\label{prop:nichols-diagonal-gkd>0} Let $\V$ be infinite-dimensional and of diagonal type with connected diagram.
If $\GK \toba(\V) < \infty$, then $\GK \toba(\V)  = 0$.
\end{prop}

We omit the proof which is analogous to the proof of Proposition \ref{prop:gk0}. Thus the classification of the braided vector spaces 
with infinite basis of diagonal type whose Nichols algebra has finite $\GK$ would be the list in Proposition \ref{prop:gk0}.

\medbreak
Our main result, Theorem \ref{theorem:main}, provides the classification of those $\V$ in the second class (braided vector spaces whose components
are blocks or points described in \S \ref{subsubsec:intro-class})
such that $\GK \toba(\V) < \infty$. This result generalizes, and is based on, \cite[Theorem 1.10]{triang}--in particular it assumes the validity of Conjecture \ref{conjecture:nichols-diagonal-finite-gkd}. In fact, the class considered here
 is an extension of that in \cite[Definition 1.8]{triang}.
As illustration, we describe examples of Nichols algebras of infinite rank
with $\GK = n$ for all $n \in \N_{\ge 2}$. 

We also observe that \cite[Theorem 1.10]{triang} does not conclude the classification of finite-dimensional braided vector spaces arising as Yetter-Drinfeld 
modules over abelian groups whose Nichols algebra has finite $\GK$, since the determination of those containing a pale block is still open, see \cite[\S 8.1]{triang}. 
Correspondingly our Theorem \ref{theorem:main} does not conclude the classification of those $\V$ as above whose Nichols algebra has finite $\GK$.

Finally, we explain how to obtain for some of these examples new pointed Hopf algebras with finite $\GK$ (albeit not finitely generated).

 \section{Preliminaries}
  
 \subsection{Conventions}\label{subsection:conventions}
  If $\ell < \theta \in\N_0$, then we set $\I_{\ell, \theta}=\{\ell, \ell +1,\dots,\theta\}$, $\I_\theta 
  = \I_{1, \theta}$. 
  Let $\G_N$ be the group of roots of unity of order $N$ in $\ku$ and $\G_N'$ the subset of primitive roots of order $N$;
  $\G_{\infty} = \bigcup_{N\in \N} \G_N$.
  All the vector spaces, algebras and tensor products  are over $\ku$.
  
  By abuse of notation, $\langle a_i: i\in I\rangle$ denotes either the group, the subgroup or the vector subspace generated by the $a_i$'s,
  the meaning being clear from the context.
  
  \smallbreak
  All Hopf algebras in this paper have bijective antipode.
  Let $H$ be a Hopf algebra.
  We refer to \cite{AS Pointed HA} for the definitions of braided vector spaces and
  the category $\ydh$ of Yetter-Drinfeld modules over $H$.
  As customary, we go back and forth between Hopf algebras in $\ydh$ and braided Hopf algebras-- that is,
  rigid braided vector spaces with compatible algebra and coalgebra structures \cite{T}.
  If $V, W \in \ydh$,  then $c_{V, W}: V\otimes W \to W \otimes V$ denotes the corresponding braiding.
  If $R$ is a Hopf algebra  in $\ydh$, then $R\# H$ is the bosonization of $R$ by $H$.

  \smallbreak We denote by $\widehat G$ the group of multiplicative characters (one-dimensional representations) of a group $G$.
  Let $\Gamma$ be an abelian group.
  The objects in $\ydG$ are the same as $\Gamma$-graded $\Gamma$-modules,
  the $\Gamma$-grading is denoted $V = \oplus_{g\in \Gamma} V_g$.
  If $g\in \Gamma$ and $\chi \in \widehat\Gamma$, then  the one-dimensional vector space $\ku_g^{\chi}$,
  with action and coaction given by $g$ and $\chi$, is in $\ydh$.
  
  Nichols algebras are graded Hopf algebras in $\ydh$, or also braided graded Hopf algebras, coradically graded and generated
  in degree one. See \cite{AS Pointed HA} for alternative characterizations.

\subsection{Convex PBW-bases and Gelfand-Kirillov dimension}\label{subsection:convex}
 
 Our reference for the notion and properties of Gelfand-Kirillov dimension is \cite{KL}.
 
 \medbreak
Let $A$ be an algebra. 
A \emph{PBW-basis} of $A$ is a $\ku$-basis $B = B(P,S,<,h)$ of $A$  that has the form
\begin{align*}
B &= \big\{p\,s_1^{e_1}\dots s_t^{e_t}:& t &\in \N_0,\ s_i \in S, \ p \in P,& s_1&>\dots >s_t, 
&  0&<e_i<h(s_i) \big\},
\end{align*}
where $P$ and $S$ are non-empty subsets of $A$; $<$ is a total order on $S$ and $h$ is a function $h: S \mapsto \N \cup \{ \infty \}$ called the height.
The elements of $S$ are called the PBW-generators. 

From now on we assume that $P = \{1\}$ and that $S$ is finite or countable with a numeration $S = \{s_1,s_2,\dots\}$ such that $i<j$ iff $s_i<s_j$.
Then we may express any $b \in B$, $b\neq 1$, as $b = s_N^{e_N}\cdots s_1^{e_1}$ where $0 \leq e_i < h(s_i)$, $i\in \I_N$, and $e_N\neq 0$; we set
\begin{align*}
\deg b &= (e_1,\dots, e_N, 0, \dots) \in \N_0^{\N}.
\end{align*}
Let $\preceq$ be the lexicographical order, reading from the right, on the set $\N_0^{(\N)}$ 
of elements of finite support of $\N_0^{\N}$ and let $\delta_{j} \in \N_0^{(\N)}$ be the element with all 0's except 1 in the place $j$. 
We consider the $\N_0^{(\N)}$-filtration on $A$ given by
\begin{align*}
A_f &= \langle s_n^{e_n}\dots s_1^{e_1} \in B: \, (e_1, e_2, \dots) \preceq (f_1, f_2, \dots) \rangle,
\end{align*}
$f = (f_1, f_2, \dots) \in \N_0^{(\N)}$. That is, $A_f = \langle b \in B: \deg b \preceq f \rangle$.

The following Definition, Lemma \ref{lema:dck} and Remark \ref{rem:GK-det} are inspired by \cite{DCK}.

\begin{definition}
The PBW-basis $B$ is \emph{convex} if  $(A_f)_{f \in \N_0^{(\N)}}$ is an algebra filtration.
\end{definition}
  
\begin{lemma}\label{lema:dck}
The PBW-basis $B$ is convex if and only if 
\begin{enumerate}[leftmargin=*,label=\rm{(\alph*)}]
\item\label{item:PBW-convex-1} for every  $i, j\in \N$ with $i <j$, there exists $\lambda_{ij} \in \ku$ such that
\begin{align}\label{eq:PBW-convex-sisj}
s_is_j &= \lambda_{ij} s_js_i + \sum_{f \prec \delta_{i} + \delta_{j}} A_f;
\end{align}
\item\label{item:PBW-convex-2} for every  $i\in \N$ such that $h(s_i)\in\N$,
\begin{align}\label{eq:PBW-convex-prv}
s_i^{h(s_i)} & \in \sum_{f \prec h(s_i) \delta_{i}} A_f.
\end{align}
\end{enumerate}
\end{lemma}

\pf
If $B$ is convex, then \ref{item:PBW-convex-1} and \ref{item:PBW-convex-2} follow directly.

Now assume that \ref{item:PBW-convex-1} and \ref{item:PBW-convex-2} hold. We claim that
\begin{align}\label{eq:claim-PBW-filtration}
s_i A_f & \subseteq A_{f+\delta_i} & \text{for all }f = (f_1, f_2, \dots) \in\N_0^{(\N)}, \, i\in\I.
\end{align} 
Let $N(f)=\max\{i\in\N: f_i\neq 0\}$. We prove the claim by induction on $N(f)$.


If $N(f)=1$, then $f=n\delta_1$ for some $n\in\N$. By \eqref{eq:PBW-convex-prv} we have that $s_1^{h(s_1)}\in A_{(h(s_1) - 1 )\delta_{1}} $ if $h(s_1)\in\N$. Hence the subalgebra generated by $s_1$ is either isomorphic to $\Bbbk[t]$ if $h(s_1)=\infty$ or else to $\Bbbk[t]/\langle m_{s_1} \rangle$ if $h(s_1)\in \N$ (where $m_{s_1}$ is the minimal polynomial of $s_1$), 
and the claim follows for $i=1$. If $i>1$, then 
$s_is_1^n\in A_{n\delta_1+\delta_i}$ by definition.

Now assume that $N:=N(f)>1$ and the claim holds for all $e$ such that $N(e)<N$. We have to prove that
\begin{align*}
s_i  s_N^{f_N} \dots s_1^{f_1} &\in A_{f+\delta_i} & \text{for all }f_i \in\N_0, \, i\in\I.
\end{align*}
Now we use induction on $f_N$. Set $f'=f-\delta_N\in\N_0^{(\N)}$. Hence $N(f')\le N$.

We assume that $f_N>1$ and the claim holds for all $e$ such that either $N(e)<N$ or else $N(e)=N$ and $e_N<f_N$.
The case $f_N=1$ follows as the recursive step. 
We have three cases. If either $i>N$ or else $i=N$ and $h(s_N)>f_N+1$, then $s_i  s_N^{f_N} s_{N-1}^{f_{N-1}} \dots s_1^{f_1} \in A_{f+\delta_i}$ by definition.

If $i=N$ and $h(s_N)=f_N+1$, then we use \eqref{eq:PBW-convex-prv} and the inductive hypothesis:
\begin{align*}
s_N^{h(s_N)} s_{N-1}^{f_{N-1}}\dots s_1^{f_1} &\in \sum_{e \prec h(s_N) \delta_{N}} A_e \, s_{N-1}^{f_{N-1}}\dots s_1^{f_1} 
\\
&= \sum_{d: N(d)< N} \sum_{j =0}^{h(s_N)-1} s_N^j A_d \, s_{N-1}^{f_{N-1}}\dots s_1^{f_1}
\\
&\subseteq \sum_{d: N(d)< N} \sum_{j =0}^{h(s_N)-1} s_N^j A_d
\subseteq  A_{h(s_N)\delta_N} \subseteq A_{f+\delta_N}.
\end{align*}

Finally, let $i<N$. By \eqref{eq:PBW-convex-sisj},
\begin{align*}
s_i s_N^{f_N} s_{N-1}^{f_{N-1}}\dots s_1^{f_1} &\in 
\lambda_{iN} s_Ns_i s_N^{f_N-1} s_{N-1}^{f_{N-1}}\dots s_1^{f_1} \\
& \qquad + \sum_{e \prec \delta_{i} + \delta_{N}} A_e s_N^{f_N-1} s_{N-1}^{f_{N-1}}\dots s_1^{f_1}.
\end{align*}
By inductive hypothesis, $s_i s_N^{f_N-1}s_{N-1}^{f_{N-1}}\dots s_1^{f_1}\in A_{f'+\delta_i}$. Thus, by definition, $s_Ns_i s_N^{f_N-1} s_{N-1}^{f_{N-1}}\dots s_1^{f_1}\in A_{f+\delta_i}$ . 

On the other hand, if $e \prec \delta_{i} + \delta_{N}$, then either $e_N=0$ or else $e_N=1$ and $e_i=\dots =e_{N-1}=0$. In the first case, by inductive hypothesis,
\begin{align*}
A_e s_N^{f_N-1} s_{N-1}^{f_{N-1}}\dots s_1^{f_1} \subseteq \sum_{d: N(d)< N} \sum_{j =0}^{f_N-1} s_N^j A_d \subseteq A_{f_N \delta_N} \subseteq A_{f+\delta_i}.
\end{align*}
In the second case, $A_e \subseteq \sum_{d: N(d)< i} s_NA_d$; by inductive hypothesis again,
\begin{align*}
A_e s_N^{f_N-1} s_{N-1}^{f_{N-1}}\dots s_1^{f_1} & 
\subseteq \sum_{d: N(d)< i} s_N A_d s_N^{f_N-1} s_{N-1}^{f_{N-1}}\dots s_1^{f_1}
\\ &
\subseteq \sum_{d: N(d)< i} s_N A_{f'+d} \subseteq s_N A_{f'+\delta_i} \subseteq A_{f+\delta_i}.
\end{align*}

Finally, from \eqref{eq:claim-PBW-filtration}, $A_eA_f\subseteq A_{e+f}$ for all $e,f\in\N_0^{(\N)}$.
\epf

\begin{remark}\label{rem:GK-det}
Assume that in \eqref{eq:PBW-convex-sisj}, $\lambda_{ij}\neq 0$ for all $i<j$. Then the associated graded algebra $\gr A$ is a (truncated) quantum linear space: $\gr A$ is the algebra presented by generators $s_i$ and relations
\begin{align*}
s_is_j &= \lambda_{ij} s_j s_i, \quad i<j, &
s_i^{h(s_i)}&=0, \quad h(s_i)<\infty.
\end{align*}
If $S$ is finite, then $\GK A = \GK \gr A= | \{s\in S: h(s)=\infty\}|$, hence $S$ is a GK-deterministic subspace of $A$, cf. \cite[Lemma 3.1]{triang}.
\end{remark}

\begin{remark}
Let $A, A'$ be subalgebras of an algebra $C$ which have convex PBW bases with PBW-generators $S$ and $S'$ respectively.
Assume that for each $s\in S$, $t\in S'$ there exists $\lambda_{s,t}\in \Bbbk$ such that $st =\lambda_{s,t} ts$, 
and that the multiplication induces a linear isomorphism $C \simeq A \otimes A'$. Then $C$ also has a convex PBW basis with PBW-generators $S \cup S'$.
\end{remark}

\begin{remark}
Let $\toba$ be a pre-Nichols algebra of a braided vector space of diagonal type.
Then $\toba$ has a convex PBW basis by \cite[Theorem 2.2]{Kh}.
Here we use the \emph{deg-lex order} \cite[\S 1.2.3]{Kh}.
\end{remark}

\begin{remark} By inspection, every Nichols algebra $\toba(V)$ with finite $\GK$ appearing in 
\cite[\S 4, 5, 7]{triang} has a convex PBW basis.
\end{remark}

\section{Locally finiteness}\label{subsec:proof-locally-finite-Z}

Recall that a family $\Fg$ of elements of a partially ordered set $(\Xx, \prec)$ is 
\emph{filtered} if given $U, W \in \Fg$, there exists $Z\in \Fg$
such that $U \prec Z$, $W \prec Z$.

For instance, the set of  Yetter-Drinfeld submodules of a given Yetter-Drinfeld module is
partially ordered by inclusion; and the  family of its finite-dimensional submodules is filtered.
A Yetter-Drinfeld module is \emph{locally finite} if it is the union of its finite-dimensional submodules.

However the family of finite-dimensional braided subspaces of a braided vector space is not necessarily filtered
(in the set of braided subspaces ordered by inclusion).

\begin{example}
Let $\triangleright: \Z \times \Z \to \Z$ be the map given by $i \triangleright j = 2i - j$, $i,j \in \Z$.
Let $V$ be a braided vector space with a basis $(x_i)_{i\in \Z}$ and braiding $c(x_i \otimes x_j) = x_{i \triangleright j} \otimes x_i$, $i,j \in \Z$.
Then the braided subspace generated by $x_0$ and $x_1$ is $V$, thus the family of finite-dimensional braided subspaces is not filtered.
\end{example}

\begin{definition}
A braided vector space is \emph{locally finite} if it is the union of its finite-dimensional braided subspaces and these form a filtered family (of the set of braided vector spaces).
\end{definition}

\begin{remark}\label{rem:loc-finite-equiv}
A braided vector space is  locally finite  if and only if every finite-dimensional subspace is contained in a  finite-dimensional braided one.
\end{remark}

If $W \in \ydh$ is a locally finite Yetter-Drinfeld module, then $W$ is a locally finite braided vector space,
but the converse is not true even if $\GK \toba (W) < \infty$.

  \begin{example}\label{example:locally-finite} Let $\Gamma = \Z^2$, $g=(1,0), h=(0,1)$.
 Let $W\in \ydG$
 with basis $(w_i)_{i\in \Z }$ such that $W=W_g$,
 $g\cdot w_i = -w_i$ and $h\cdot w_i=w_{i+1}$ for all $i\in \Z$.
 Then the action of $\Gamma $ on $W$ is not locally finite, but $W$  is a locally finite braided vector space,
 $\toba (W)=\Lambda W$ and $\GK \toba (W)=0$. Furthermore, $\toba (W) \# \ku \Gamma = \ku \langle g^{\pm 1}, h^{\pm 1},w_0 \rangle$ is a finitely generated graded algebra. By \cite[Example 2.4]{triang} (compare with \cite[Lemma 2.2]{triang}),
 \begin{align*}
 \GK (\toba (W) \# \ku \Gamma)  = \infty > 2 = \GK \toba (W) + \GK \ku \Gamma.
 \end{align*}
  \end{example}

  \begin{question}\label{question:nichols-locally-finite}
 Let $(V, c)$ be a braided vector space with $\GK \toba(V) < \infty$.
 Is $(V, c)$ a locally finite braided vector space?
  \end{question}
  
  The defining relations of Nichols algebras of locally finite braided vector spaces
  could be determined from their finite-dimensional counterparts by the following  fact.
  If $(V,c)$ is a braided vector space, then set $\cJ(V) =$ the ideal of $T(V)$ of relations in $\toba(V)$.

  \begin{lemma}\label{prop:infinite-dynkin-rels}
 Let $(V, c)$ be a  braided vector space and let $\Fg$ be a filtered family of braided subspaces such that
 $V = \displaystyle\bigcup_{W\in \Fg} W$. Then
 \begin{align}\label{eq:infinite-dynkin-rels}
 \cJ(V) &= \displaystyle\bigcup_{W\in \Fg} \cJ(W), 
 \\ \label{eq:infinite-dynkin-toba}
 \toba(V) &= \displaystyle\bigcup_{W\in \Fg} \toba(W), 
 \\ \label{eq:infinite-dynkin-toba-gkdim}
 \GK \toba(V) &= \displaystyle\sup_{W\in \Fg} \GK \toba(W).  
 \end{align}
  \end{lemma}
  
  \pf Since $\Fg$ is filtered, the right-hand side $\cI$ of \eqref{eq:infinite-dynkin-rels} 
  is a homogeneous Hopf ideal of $T(V)$
  that intersects $\ku \oplus V$ trivially, hence $\cI \subseteq \cJ(V)$.
  Conversely, let $r\in \cJ(V)$. Then there exists $W\in \Fg$ such that $r\in T(W)$.
  Now $\cJ(V) \cap T(W)$ is a homogenous Hopf ideal of $T(W)$ that intersects $\ku \oplus W$ trivially,
  hence $\cJ(V) \cap T(W) \subset \cJ(W)$ and $r\in \cJ(W)$. Thus $\cI \supseteq \cJ(V)$.
  The contention $\subseteq$ in \eqref{eq:infinite-dynkin-toba} is immediate as $\toba(V)$ is generated by $V$,
  and we also know that $\toba(W) \hookrightarrow \toba(V)$ for any braided subspace $W$. Finally 
  \eqref{eq:infinite-dynkin-toba} implies \eqref{eq:infinite-dynkin-toba-gkdim} using again that $\Fg$ is filtered.
  \epf
  
  \begin{coro}\label{cor:infinite-dynkin-rels}
 Let $(V, c)$ be a locally finite braided vector space. Then
 \begin{align}\label{eq:infinite-dynkin-rels-cor}
 \cJ(V) &= \displaystyle\bigcup_{\substack{W\ \text{\rm braided subspace}\\ \dim W < \infty}} \cJ(W), 
 \\ \label{eq:infinite-dynkin-toba-gkdim-cor}
 \GK \toba(V) &= \displaystyle\sup_{\substack{W\ \text{\rm braided subspace}\\ \dim W < \infty}} \GK \toba(W).
 \end{align}
  \end{coro}
  
  \pf Apply Lemma \ref{prop:infinite-dynkin-rels} to  the family of finite-dimensional braided subspaces  of $V$.
  \epf 
  
\subsection{Local finiteness over an abelian group}
Let now $\Gamma$ and $\V$ be as in the Introduction. Then $\V = \oplus_{g\in \Gamma} \V_g \in \ydG$.

\begin{theorem}\label{prop:locally-finite-Z}
If $\GK \toba(\V) < \infty$, then $\V$ is a locally finite braided vector space.
\end{theorem}

\pf If $v = \sum v_g\in \V$, where $v_g \in \V_g$, then  $\supp v := \{g \in \Gamma: v_g \neq 0 \}$.
Also, $\supp \V = \{g \in \Gamma: \V_g \neq 0 \}$.

\begin{paso}\label{step:h-not-loc-finite}
Let $V$ be a $\Gamma$-module. If the action of $h\in \Gamma$ on $V$ is not locally finite,
then there is $v \in V$ such that
$(v_n)_{n\in \Z}$ is linearly independent, $v_n := h^n \cdot v$.
\end{paso}

Let $v \in V$ and $v_n := h^n \cdot v$. 
Assume that there is a non-trivial relation $\sum_{q\le n\le p} a_n v_n = 0$
with $p-q$ minimal; applying $h^{-q}$, we may assume that $q = 0$. Then $\langle v_n: 0\le n\le p-1 \rangle$
is both stable by the action of $h$ and $h^{-1}$. If this happens for all $v\in V$, then the action of $h$
is locally finite.

\begin{paso}\label{step:g-Vg-loc-finite}
Let $g\in \supp \V$. Then the action of $g$ on $\V_g$ is locally finite.
\end{paso}

Otherwise, by Step \ref{step:h-not-loc-finite}, there is $v \in \V_g$
such that $(v_n)_{n\in \Z}$ is linearly independent, $v_n := g^n \cdot v$. Now 
$U:= \langle v_n:n\in \Z \rangle$ is a braided vector subspace of $\V$, since $c(v_n \otimes v_m) = v_{m + 1} \otimes v_n$. 
Then $\toba(U) = T(U)$. Indeed, let $\Omega_n = \sum_{\sigma \in \Sb_n} M_{\sigma}\in \End T^n(U)$ be the quantum symmetrizer,
so that $\toba^n(U) = T^n(U) / \ker\Omega_n$, and consider the $\Z$-grading in $T(U)$ given by $U_m = \langle v_m\rangle$, $m \in \Z$.
Then $M_{\sigma} (T^n(U)_m) = T^n(U)_{m + \ell(\sigma)}$ for all $m\in \Z$, where $\ell$ is the usual length. Since $\Sb_n$ has a unique element
of maximal length, $\Omega_n$ is injective, proving our claim. Whence $\toba(U)$ and \emph{a fortiori} $\toba(\V)$ have infinite $\GK$.

\bigbreak 
Let $g\in \supp \V$ and $\lambda\in \kut$. As usual we set $\V_g^{\lambda} =\{v\in \V_g: g\cdot v= \lambda v \}$,
$\V_g^{(\lambda)} =\{v\in \V_g: (g -\lambda)^n \cdot v= 0 \text{ for } n >>0\}$. 
By Step \ref{step:h-not-loc-finite}, $\V_g  = \oplus_{\lambda\in \kut} \V_g^{(\lambda)}$.
Since $\Gamma$ is abelian, its action preserves $\V_g^{\lambda}$ and $\V_g^{(\lambda)}$ for all $\lambda$.

\begin{paso}\label{step:h-Vglambda-loc-finite}
Let $g, h\in \supp \V$ and $\lambda\in \kut$. Then the action of $h$ on $V_g^{\lambda}$ is locally finite.
\end{paso}

Otherwise, by Step \ref{step:h-not-loc-finite}, there is $v \in \V_g^{\lambda}$
such that $(x_n)_{n\in \Z}$ is linearly independent; here $x_n := h^n \cdot v$. 
Now $c(x_n \otimes x_m) = \lambda x_{m} \otimes x_n$ so that $\lambda = -1$ by 
\cite[Proposition 3.1]{diag}. We distinguish two cases:

\begin{enumerate}[leftmargin=*,label=\rm{(\Alph*)}]
\item\label{item:g-Vh-locfin} The action of $g$ on $\V_h$ is locally finite.

\item\label{item:g-Vh-not-locfin} The action of $g$ on $\V_h$ is not locally finite.
\end{enumerate}

Assume \ref{item:g-Vh-locfin}. By Step  \ref{step:g-Vg-loc-finite}, there are $y \in \V_h$
and $\mu, q \in \kut$ such that
\begin{align*}
g\cdot y &= \mu y, & h\cdot y &= q y.
\end{align*}
We may assume that $\V = \V_g \oplus \V_h$, $\V_g = \langle x_i: i\in \Z \rangle$ and $\V_h = \langle y \rangle$. 
Let $n\in\N_{\ge 2}$. Then $K[n] = \langle x_r - x_{r+n}: r\in \Z \rangle$
is a Yetter-Drinfeld submodule of $\V$ and $\Uu[n] = \V / K[n] =  \Uu[n]_g \oplus \Uu[n]_h$,
where $\{\overline{x_i}\}_{i\in\I_n}$ is a basis of $\Uu[n]_g$ and $\Uu[n]_h = \V_h$. We fix $\zeta\in\G_n'$ and set $z_j=\sum_{i \in \I_n} \zeta^{-ij} \overline{x_i}$. Thus
\begin{align*}
g\cdot z_j &= -z_j, & 
h\cdot z_j &= \sum_{i \in \I_n} \zeta^{-ij} \overline{x_{i+1}} = \zeta^j \, z_j,
\end{align*}
and the braiding of $\Uu[n]$ satisfies
\begin{align*}
c(z_i\otimes z_j) &= -z_j\ot z_i, &
c(z_i\otimes y) &=  \mu \, y \ot z_i, \\
c(y\otimes z_j) &= \zeta^j \, z_j\ot y, &
c(y\otimes y) &= q \, y\ot y.
\end{align*}
That is, $\Uu[n]$ is of diagonal type with diagram
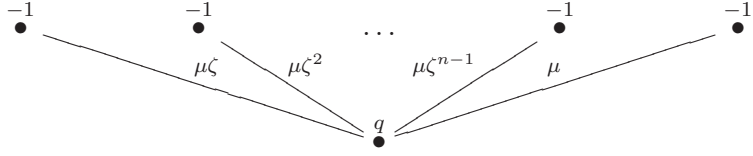
\begin{align*}
\xymatrix@C=50pt{
\overset{-1}{\bullet}  & \overset{-1}{\bullet} & \dots &
\overset{-1}{\bullet}  & \overset{-1}{\bullet}   
\\
& & \overset{q}{\bullet} \ar@{-}[ull]_{\mu\zeta} 
\ar@{-}[ul]_{\mu\zeta^2}  \ar@{-}[ur]^{\mu\zeta^{n-1}} 
\ar@{-}[urr]^{\mu} & & }
\end{align*}

We consider five cases:
\begin{itemize}[leftmargin=*]
\item $q\notin \G_{\infty}$. Suppose that $\GK \toba(\Uu[n])<\infty$. By \cite[Lemma 3.3]{diag} there exist $h,j\in\N_0$ such that $q^h(\mu\zeta) = 1 = q^j\mu$, so $q^{j-h}\in\G_{n}'$, a contradiction.
\item $q\in \G_{N}'$, $N>24$. Suppose that $\GK \toba(\Uu[2N])<\infty$. By \cite[Theorem 4.1]{diag}, each subdiagram of rank two appears in \cite[Table 1]{H-classif}. Thus, for each $i\in\I_{2N}$ there exists $a_{i}\in\{-3,-2,-1,0\}$ such that $q^{a_i} (\mu\zeta^i)=1$. Hence $\zeta=q^{a_1-a_2}\in \G_N$, a contradiction.
\item $q\in \G_{N}'$, $3\le N \le 24$. Suppose that $\GK \toba(\Uu[100])<\infty$. By \cite[Theorem 4.1]{diag}, each subdiagram of rank two appears in \cite[Table 1]{H-classif}: By inspection the labels of the edges belong to $\bigcup_{k\in\I_{3,72}} \G_k$. Thus $\zeta=(\mu\zeta)\mu^{-1}\in\bigcup_{k\in\I_{3,72}} \G_k$, a contradiction.
\item $q=-1$. Suppose that $\GK \toba(\Uu[n])<\infty$ for all $n\ge 2$. 
If $\mu\neq \pm 1$, then the braided vector space $\cR_n(\Uu[n])$ obtained by reflection at the vertex $x_n$ has diagram
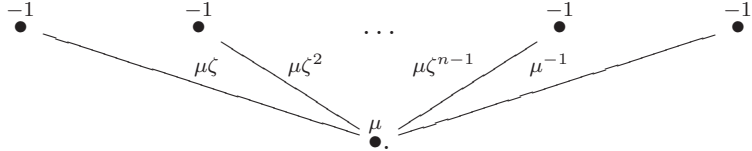
\begin{align*}
\xymatrix@C=50pt{
\overset{-1}{\bullet}  & \overset{-1}{\bullet} & \dots &
\overset{-1}{\bullet}  & \overset{-1}{\bullet}   
\\
& & \overset{\mu}{\bullet}. \ar@{-}[ull]_{\mu\zeta} 
\ar@{-}[ul]_{\mu\zeta^2}  \ar@{-}[ur]^{\mu\zeta^{n-1}} 
\ar@{-}[urr]^{\mu^{-1}} & & }
\end{align*}
A similar work as in the previous cases shows that $\GK \toba(\cR_n(\Uu[n]))=\infty$ for some $n$, depending on the case. 
When $\mu=\pm 1$, we consider $\cR_1(\Uu[n])$ and conclude that $\GK \toba(\cR_1(\Uu[n]))=\infty$ by an analogous analysis.
But this is a contradiction with \cite[Theorem 2.4]{diag}.
\item $q=1$. Here $\GK \toba(\Uu[2])=\infty$, since either $\mu\neq 1$ or else $-\mu\neq 1$, and then \cite[Lemma 2.8]{triang} applies for a subspace of $\Uu[2]$.
\end{itemize}

In any case there exists $n\geq 2$ such that $\GK \toba(\Uu[n])=\infty$, hence $\GK \toba(\V)=\infty$.

\bigbreak
Assume \ref{item:g-Vh-not-locfin}.
Then, by Step \ref{step:h-not-loc-finite}, there is $w \in \V_h^{q}$
such that $(y_r)_{r\in \Z}$ is linearly independent; here $y_r := g^r \cdot w$.
Thus we may assume that $\V = \V_g \oplus \V_h$, $\V_g = \langle x_i: i\in \Z \rangle$ 
and $\V_h = \langle y_r: r\in \Z \rangle$. Now $K = \langle y_r - y_{r+1}: r\in \Z \rangle$
is a Yetter-Drinfeld submodule of $\V$ and $\Uu = \V / K =  \Uu_g \oplus \Uu_h$,
where $\Uu_g = \V_g$ and $\dim \Uu_h = 1$. So, we are in the situation \ref{item:g-Vh-locfin}.

\begin{lemma}\label{lemma:ext-locfin} 
Let $\varUpsilon = \langle \gamma_1, \dots, \gamma_r\rangle$ be a finitely generated abelian group 
and  let
$\xymatrix{0 \ar  @{->}[r] & W'\ar  @{->}[r]& W \ar  @{->}[r]& W'' \ar  @{->}[r] & 0}$ be an exact sequence of $\varUpsilon$-modules.
If the actions of $\varUpsilon$ on $W'$ and $W''$ are locally finite, then so is the action of $\varUpsilon$ on $W$.
\end{lemma}

\pf Let $w \mapsto \overline w$ denote the projection $W \to W''$. Pick $w\in W$. Then
there is a $\varUpsilon$-stable submodule $U = \langle \overline{w}_1, \dots, \overline{w}_{\ell} \rangle$ of $W''$
such that $\overline w \in U$. That is, there are scalars $\alpha_i, \beta^{i}_{kj}\in \ku$ such that
\begin{align*}
\overline w &= \sum_i \alpha_i \overline{w_i}, &
\gamma_k \cdot \overline w_j &= \sum_i \beta^{i}_{kj} \overline{w_i}
\end{align*} 

Hence there are $v_0, v_{kj} \in W'$ such that
\begin{align*}
w &= \sum_i \alpha_i w_i + v_0, & \gamma_k \cdot  w_j &= \sum_i \beta^{i}_{kj} w_i +  v_{kj},& j&\in \I_{\ell}, \, k\in \I_{r}.
\end{align*}

Let $Z$ be a finite-dimensional $\varUpsilon$-submodule of $W'$ containing
$v_0$, and all the $v_{kj}$'s. Then $Z + \langle (w_i)_{i\in \I_{\ell}} \rangle$ is a finite-dimensional $\varUpsilon$-submodule of $W$ that contains $w$.  \epf

From now on, we assume without loss of generality that $\Gamma$ is generated by $\supp \V$.

\begin{paso}\label{step:fingen-loc-finite}
The action of any finitely generated subgroup of $\Gamma$ on $\V$ is locally finite.
\end{paso}

Let $\varUpsilon = \langle h_1, \dots, h_r\rangle$ be a finitely generated subgroup of $\Gamma$; we may assume that $h_1, \dots, h_r \in \supp V$.
We first claim that the action of  $\varUpsilon$ on $\V_g$ is locally finite for any $g\in \Gamma$.
Indeed, by Zorn there is a maximal locally finite $\varUpsilon$-submodule $W'$ of $\V_g$. If $W' \neq \V_g$, then consider $\Uu_g = \V_g / W'$,
$\Uu_h = \V_h$ for $h \in \supp \V \cap \varUpsilon$ and $\Uu = \Uu_g \oplus \oplus_{h\in \varUpsilon}\Uu_h$. 
By induction on the number $r$ of generators and using Step \ref{step:h-Vglambda-loc-finite}, 
there exists $\chi \in \widehat{\varUpsilon}$ such that $\Uu_g^{\chi} \neq 0$.
Pick $w \in \Uu_g^{\chi} - 0$ and set $W'' = \ku w$. Let $W$ be the submodule of $\V_g$ generated by $W'$ and a pre-image of $w$.
Then $W$ is a locally finite $\varUpsilon$-submodule of $\V_g$, contradicting the maximality of $W'$, by Lemma \ref{lemma:ext-locfin}.
This shows the claim and a standard argument gives the Step.

\begin{paso}\label{step:fingen-bvs}
$\V$ is a locally finite braided vector space.
\end{paso}

By Remark \ref{rem:loc-finite-equiv}, it is enough to consider a vector subspace $\mathcal V = \langle v_1, \dots, v_m \rangle$. Let 
$S = \{h_1, \dots, h_r\} = \bigcup_{i\in \I_m} \supp v_i$, $\varUpsilon = \langle h_1, \dots, h_r\rangle$ and 
$\V_S = \oplus_{h \in S} \V_h$. Then $\V_S$ is a locally finite Yetter-Drinfeld module over $\ku \Upsilon$ by Step \ref{step:fingen-loc-finite}, 
hence it is a locally finite braided vector space, and $\mathcal V$ is a contained in a finite-dimensional braided vector space.
\epf

  \subsection{Decompositions}
  As in \cite[Definition 2.1]{Grana}, a decomposition of a braided vector space
  $V$ is a family $(V_i)_{i\in I}$ of subspaces such that
  \begin{align}\label{eq:bradinig-generalform}
  V &=  \oplus_{i \in I} V_i,
  & c(V_i \otimes V_j) &=  V_j \otimes V_i,&  i,j&\in I.
  \end{align}
  If $i \neq j \in I$, then we set
  $c_{ij} = c_{\vert V_i \otimes V_j}: V_i \otimes V_j \to V_j \otimes V_i$.
  
\begin{question}\label{question:cijcji-id}
Assume that $c_{ji}c_{ij} = \id_{V_i \otimes V_j}$ for all $i \neq j \in I$. Is it true that
    \begin{align}\label{eq:cijcji-id}
    \GK \toba(V)& =  \sum_{i \in I} \GK\toba(V_i)?
    \end{align}
If yes, then $\GK \toba(V) < \infty$ implies $\GK \toba(V_i)  = 0$ for all but finitely many  $i\in I$.
\end{question}

  If $F \subset I$, then $V_F =  \oplus_{i \in F} V_i$ is a braided subspace of $V$.
  Hence 
  \begin{align*}
  \Fg &=  \left\{V_F: F \subset I, \, \vert F \vert < \infty \right\}
  \end{align*}
  is a filtered family of braided subspaces of $V$; by Lemma \ref{prop:infinite-dynkin-rels}, we have
  \begin{align}\label{eq:infinite-dynkin-rels-decomposition}
  \cJ(V) &= \displaystyle\bigcup_{F \subset I, \, \vert F \vert < \infty} \cJ(V_F), 
  \\ \label{eq:infinite-dynkin-toba-gkdim-decomposition}
  \GK \toba(V) &= \displaystyle\sup_{F \subset I, \, \vert F \vert < \infty} \GK \toba(V_F).
  \end{align}
  
  Let $F \subset I$, $\vert F \vert < \infty$; fix an ordering $i_1, \dots, i_k$ of $F$. 
  By the proof of \cite[2.2]{Grana}, the multiplication induces a monomorphism of graded vector spaces
  \begin{align}\label{eq:finite-decomposition-inclusion}
  \toba(V_{i_1}) \otimes \toba(V_{i_2}) \otimes \dots \otimes \toba(V_{i_k}) & \hookrightarrow \toba(V).
  \end{align}

\begin{question}\label{question:cijcji-gral}
Is it true that
  \begin{align}\label{eq:cijcji-gral}
  \GK \toba(V)& \geq  \sum_{i \in F} \GK\toba(V_i)?
  \end{align}
Assuming that $\dim V_i < \infty$ for all $i\in F$? Assuming this and that the Hilbert series of  
$\toba(V_i)$ is rational for all $i\in F$?
\end{question}

  \section{Diagonal type}\label{section:diagonal-type}
  A \emph{point} of label $q\in \kut$ is a braided vector space $(V,c)$ of dimension 1 with $c = q\id$.
  Let $V$ be a braided vector space of diagonal type; that is there are $(x_i)_{i\in I}$ a basis of $V$ and 
  $\bq = (q_{ij})_{i,j\in I}\in \ku^{I \times I}$
  such that $q_{ij}\ne 0$ and $c(x_i\otimes x_j)=q_{ij}x_j\otimes x_i$ for all $i,j\in I$.
    It turns out that there are interesting examples of infinite-dimensional
  braided vector spaces of diagonal type with $\GK = 0$, so we also ask:
  
  \begin{question}\label{question:finite-GK-diagonal-infinitematrix}
 Classify all braided vector spaces $(V, c)$ of diagonal type
 with matrix  $(q_{ij})_{i,j\in I}$, where $I$ is infinite countable  such that
 $\GK \toba(V) < \infty$.
  \end{question}

We first describe two classes of infinite-dimensional braided vector spaces $(V,c)$ of diagonal type.
Let $I = \N =\{1, 2, \dots\}$ or $\Z$.
First, consider $\ba = (a_{ij})_{i,j\in I}$ with Dynkin diagram as in Table \ref{tab:dynkin-infinite}.
Then $(q_{ij})_{i,j\in I}$ is of Cartan type $\ba$ if $q_{ii}\neq  1$ and $q_{ij}q_{ji} = q_{ii}^{a_{ij}}$ holds for all $i\neq j \in \I_{\theta}$.
  
  \begin{table}[ht]
 \caption{Infinite Dynkin diagrams}\label{tab:dynkin-infinite}
 \begin{center}
  \begin{tabular}{c l}
  $A_{\infty}$ & $\xymatrix@C-15pt{\cdots \ar  @{-}[r]  & \circ\ar  @{-}[r]  &
   \circ\ar@{.}[r]&  \circ\ar  @{-}[r] & \circ \ar  @{-}[r]  & \circ \ar  @{-}[r]  & \cdots}$
  \\
  $A_{+\infty}$ & $\xymatrix@C-15pt{  \circ\ar  @{-}[r]  &
   \circ\ar@{.}[r]&  \circ\ar  @{-}[r] & \circ \ar  @{-}[r]  & \circ \ar  @{-}[r]  & \cdots}$
  \\
  $B_{\infty}$ & $\xymatrix@C-15pt{  \circ \ar  @{<=}[r]  & \circ \ar  @{-}[r]  &
   \circ\ar@{.}[r]&  \circ\ar  @{-}[r] & \circ \ar  @{-}[r]  & \cdots}$
  \\
  $C_{\infty}$ & $\xymatrix@C-15pt{  \circ  \ar  @{=>}[r]  & \circ\ar  @{-}[r]  &
   \circ\ar@{.}[r]&  \circ\ar  @{-}[r] & \circ \ar  @{-}[r]  & \cdots}$
  \\
  \raisebox{-20pt}{$D_{\infty}$} & $\xymatrix@R-15pt@C-15pt{  & \circ \ar  @{-}[d] & &  & \\
   \circ\ar  @{-}[r]  &  \circ\ar@{.}[r]&  \circ\ar  @{-}[r] & \circ \ar  @{-}[r]  & \circ\ar  @{-}[r]  & \cdots}$
  \end{tabular}
 \end{center}
  \end{table}
  
  To describe the second class, we recall  from \cite{H-classif}
  that the generalized Dynkin diagram of a matrix  $(q_{ij})_{i,j\in I}$
  such that $q_{ii} \neq 1$ is a graph with set of points $I$, with the following decoration:
  
  \begin{itemize}[leftmargin=*]\renewcommand{\labelitemi}{$\circ$}
 \item The vertex $i$ is decorated with $q_{ii}$ above; when the numeration of the vertex is needed, it is stated below.
 
 \item Let $i\neq j \in I$. If $q_{ij}q_{ji} = 1$, there is no edge between $i$ and $j$, otherwise there is an edge decorated with $q_{ij}q_{ji}$.
  \end{itemize}
  
  All the diagrams in Table \ref{tab:super-infinite} obey the following conventions:
  
  \begin{itemize}[leftmargin=*]\renewcommand{\labelitemi}{$\diamond$}
 \item $q\in \ku$, $q\neq \pm 1$.
 
 \medbreak
 \item $\bp: I \to \G_2$ is a function, $\bp\not \equiv 1$; if $I = \N$, then $d = \min \{i\in I: \bp = -1\}$.
 Then $q_{ii} = \begin{cases}
 -1, &\text{if } \bp_i = -1, \\ q \text{ or } q^{-1}, &\text{if } \bp_i = 1,\end{cases}$ (unless explicitly stated).
 
 \medbreak
 \item They are locally  of the following forms (unless explicitly stated)
 \begin{align*}
 &\xymatrix@R-18pt{ \ar@{-}[r]^{q^{-1}} & \overset{q}{\circ} \ar@{-}[r]^{q^{-1}} &\hspace{-3pt},}&
 &\xymatrix@R-18pt{ \ar@{-}[r]^{q} & \overset{q^{-1}}{\circ} \ar@{-}[r]^{q} &\hspace{-3pt},} &
 &\xymatrix@R-18pt{ \ar@{-}[r]^{q} & \overset{-1}{\circ} \ar@{-}[r]^{q^{- 1}} &\hspace{-3pt},} &
 &\xymatrix@R-18pt{ \ar@{-}[r]^{q^{-1}} & \overset{-1}{\circ} \ar@{-}[r]^{q} & \hspace{-3pt}.}
 \end{align*}
  \end{itemize}
  
  The following matrices $(q_{ij})_{i, j\in I}$ give rise to braided vector spaces $(V,c)$
  with $\GK\toba(V) = 0$, being unions of finite-dimensional Nichols algebras:
  
  \begin{enumerate}\renewcommand{\theenumi}{\alph{enumi}}
 \renewcommand{\labelenumi}{(\theenumi)}
 \item
 $(q_{ij})_{i, j\in I}$  of Cartan type as in  Table \ref{tab:dynkin-infinite},
 and $q \in \G_{\infty} - 1$ for all $i\in I$.
 
 \item   $(q_{ij})_{i, j\in I}$  of super type
 as in Table \ref{tab:super-infinite},   and $q \in \G_{\infty} - 1$ for all $i\in I$.
  \end{enumerate}

  \begin{table}
 \caption{Some generalized Dynkin diagrams}\label{tab:super-infinite}
 \begin{center}
  \begin{tabular}{c c c c}
  Label & $I$ & Diagram & Parameter
  \\
  $A_{+\infty}(\bp, q)$ & $\N$ & $\xymatrix@C-10pt{ {\underset 1\circ}\ar  @{-}[r]  &
   \circ\ar@{.}[r]&  \circ\ar  @{-}[r] &\overset{-1}{\underset d\circ}\ar  @{-}[r]^{q} \ar  @{-}[r]  & \circ \ar  @{-}[r]  & \cdots}$ &
  \\
  $A_{\infty}(\bp, q)$ & $\Z$ & $\xymatrix@C-10pt{&\cdots \ar  @{-}[r]  & \circ\ar  @{-}[r]  &
   \circ\ar@{-}[r]&   \overset{-1}{\underset 0\circ}\ar  @{-}[r]^{q} & \circ \ar  @{-}[r]    & \cdots}$ &
  \\
  $B_{\infty}(\bp, \nu)$& $\N$ & $\xymatrix@C-10pt{  \overset{\nu}\circ \ar  @{-}[r]^{q^{-1}}  & \circ \ar  @{-}[r]  &
   \circ\ar@{.}[r]&  \circ\ar  @{-}[r] & \circ \ar  @{-}[r]  & \cdots}$ & $q = \nu^2$
  \\
  $B_{\infty}(\bp, \zeta)$& $\N$ & $\xymatrix@C-10pt{\overset{\zeta}{\circ} \ar  @{-}[r]^{-\zeta}  & \circ \ar  @{-}[r]  &
   \circ\ar@{.}[r]&  \circ\ar  @{-}[r] & \circ \ar  @{-}[r]  & \cdots}$ & $\zeta\in \G'_3$, $q = -\zeta^2$
  \\
  $C_{\infty}(\bp, q)$ & $\N$ & $\xymatrix@C-10pt{  \overset{q^2}\circ \ar  @{-}[r]^{q^{-2}}  & \overset{q}{\circ} \ar  @{-}[r]  &
   \circ\ar@{.}[r]&  \circ\ar  @{-}[r] & \circ \ar  @{-}[r]  & \cdots}$ & $q^4 \neq 1$
  \\
  \raisebox{-20pt}{$D_{\infty}(\bp, q)$} & \raisebox{-20pt}{$\N$}
  & $\xymatrix@R-15pt@C-10pt{  & \overset{q}{\circ} \ar  @{-}[d]^{q^{-1}} & &  & \\
   \overset{q}{\circ} \ar  @{-}[r]^{q^{-1}}  &  \circ\ar@{.}[r]&  \circ\ar  @{-}[r] & \circ \ar  @{-}[r]  & \circ\ar  @{-}[r]  & \cdots}$ &\\
  \raisebox{-20pt}{$D_{\infty}(\bp, q)$} & \raisebox{-20pt}{$\N$}
  & $\xymatrix@R-15pt@C-10pt{  & \overset{-1}{\circ} \ar  @{-}[d]^{q^{-1}}\ar  @{-}[dl]_{q^{2}} & &  & \\
   \overset{-1}{\circ} \ar  @{-}[r]_{q^{-1}}  &  \circ\ar@{.}[r]&  \circ\ar  @{-}[r] & \circ \ar  @{-}[r]  & \circ\ar  @{-}[r]  & \cdots}$ &
  \end{tabular}
 \end{center}
  \end{table}

  \bigbreak
  Conversely, let $V$ be a braided vector space of diagonal type with connected braiding, with a basis
  $(x_i)_{i\in I}$ such that $c(x_i\otimes x_j)=q_{ij}x_j\otimes x_i$,
  where  $q_{ij} \in \ku^{\times}$ for all $i,j\in I$.

  \begin{prop}\label{prop:gk0}  If $\GK \toba(V) = 0$, then either of the following holds:
 
 \begin{enumerate}\renewcommand{\theenumi}{\alph{enumi}}
  \renewcommand{\labelenumi}{(\theenumi)}
  
  \item\label{prop:gk0-infinite-diag-cartan} $(q_{ij})_{i, j\in I}$ is of Cartan type
  $A_{\infty}$ ($I=\Z$), $A_{+\infty}$, $B_{\infty}$, $C_{\infty}$, or $D_{\infty}$ ($I = \N$), see Table \ref{tab:dynkin-infinite},
  and $q_{ii} \in \G_{\infty} - 1$ for all $i\in I$.

  \item\label{prop:gk0-infinite-diag-super}   $(q_{ij})_{i, j\in I}$ is of super type
  $A_{+\infty}(\bp, q)$, $A_{\infty}(\bp, q)$, $B_{\infty}(\bp, q)$, $C_{\infty}(\bp, q)$,
  or $D_{\infty}(\bp, q)$, see Table \ref{tab:super-infinite}, and $q \in \G_{\infty} - 1$ for all $i\in I$.
 \end{enumerate}
  \end{prop}

  \pf
  The argument is standard \cite[Ex. 4.14]{K}.
  Suppose that $V$ is of Cartan type. If the Dynkin diagram contains a point $P$ with three concurrent edges, then
  $V$ contains a connected braided subspace $U$ of dim $m \ge 7$ whose diagram contains $P$; hence $U$ is of Cartan
  type $D_{m}$ by \cite{H-classif}. Now since $V$ has a connected braiding, one constructs recursively a
  braided subspace $W$ of Cartan type $D_{\infty}$. If $W \neq V$, then there is a point out of $W$ connected to a point
  in $D_{\infty}$, but this contradicts  \cite{H-classif}. So, $V = W$ is of Cartan type $D_{\infty}$.
  The argument in all other cases is analogous. \epf

  \begin{example}
 $B_{\infty}(\bp, \zeta)$, $\zeta\in \G'_3$, $q = -\zeta^2$. A set of defining relations of $\toba(V)$ 
is the union of those of the various $B_{\theta}(\bp, \zeta)$ described in \cite[6.1.4]{AA-diag-survey}, see also \cite{Ang-crelle}, 
using Lemma \ref{prop:infinite-dynkin-rels}. 
  \end{example}

  \section{Decompositions whose components are blocks or points}
  
  \subsection{Blocks}
  Let $\epsilon\in \kut$ and $\ell \in \N_{\ge 2}$. 
  Let $\cV(\epsilon,\ell)$ be the braided vector space
  with a basis $(x_i)_{i\in\I_\ell}$ such that
  \begin{align*}
  c(x_i \ot  x_1) &= \epsilon x_1 \ot  x_i,& c(x_i \ot  x_j) &=(\epsilon x_j+x_{j-1}) \ot  x_i,& i &\in \I_\ell, \, j \in \I_{2,\ell}.
  \end{align*}
  This  braided vector space is a called a \emph{block}.

  \begin{theorem}\label{theorem:blocks} \cite[Theorem 1.2]{triang}
 $\GK \toba(\cV(\epsilon,\ell)) < \infty$ if and only if $\ell = 2$ and $\epsilon \in \{\pm 1\}$,
 in which case $\GK \toba(\cV(\epsilon,\ell)) = 2$.
  \end{theorem}
  
  \subsection{A class of braided vector spaces}\label{subsubsec:intro-class}
  We consider in this Subsection braided vector spaces $(V, c)$ of the following sort. 
  Let $I$ be an infinite subset of $\Q$ such that $I \cap (I + \frac{1}{2}) = \emptyset$. We suppose that
  
  \begin{enumerate}[leftmargin=*,label=\rm{(\Alph*)}]
 \item\label{item:decomposition} $V$ has a decomposition $  V =  \oplus_{i \in I} V_i$ as in 
 \eqref{eq:bradinig-generalform}. Furthermore, 
 there exists $\emptyset \neq J \subseteq I$ such that
 $V_j \simeq \cV(\epsilon_j,\ell_j)$ is a block,  $j\in J$.
 Also, if $i\in I - J$, then $V_{i}$ is a $q_{ii}$-point, with $q_{ii} \in \ku^{\times}$; we  
 fix $x_i\in V_i - 0$, $i \in I - J$.
  \end{enumerate}
  
  Let $J_{\pm} = \{j\in J:  \epsilon_j = \pm 1,\, \ell_j = 2\}$.
  By Theorem \ref{theorem:blocks}, we may (and will) assume that $J = J_{+} \cup J_{-}$.
  Given $j \in J$, we fix a basis $B_j = \{x_j, x_{\fjdos}\}$ of $V_j$
  such that the braiding is given  by
  \begin{align*}
  (c(x_r \otimes x_s))_{r, s\in B_{j}} &=
  \begin{pmatrix}
  \epsilon_j x_j  \otimes x_j &  (\epsilon_j x_{\fjdos} + x_j ) \otimes x_j 
  \\
  \epsilon_j x_j  \otimes x_{\fjdos} & (\epsilon_j x_{\fjdos} + x_j ) \otimes x_{\fjdos}  
  \end{pmatrix}.
  \end{align*}
  
  If $i,h \in I - J$, then the braiding $c_{ih}$ is uniquely determined by $q_{ih} \in \ku^{\times}$: 
  $c_{ih} = q_{ih} \tau$, where $\tau$ is the usual flip. 
  Let
  \begin{align*}
  V_{\diag} = \oplus_{i \in I - J} V_i.
  \end{align*}

  Our next assumption deals with the braidings between blocks and points.
  
  \begin{enumerate}[leftmargin=*,label=\rm{(\Alph*)}]\addtocounter{enumi}{1}
 \item\label{item:braiding-points-blocks}
 For every  $j\in J$ and $i \in I - J$, there exist
 $q_{ij}, q_{ji} \in \kut$ and $a_{ij} \in \ku$ such that the braiding between $V_j$ and $V_i$ is given  by
 \begin{align}\label{eq:braiding-point-block}
 c(x_j \otimes x_i) &= q_{ji} x_i  \otimes x_j, & c(x_{\fjdos} \otimes x_i) &= q_{ji} x_i  \otimes x_{\fjdos}, 
 \\ \label{eq:braiding-block-point}
 c(x_i \otimes x_j) &=q_{ij} x_j  \otimes x_i ,& c(x_i \otimes x_{\fjdos}) &=  q_{ij}(x_{\fjdos} + a_{ij} x_j ) \otimes x_i.
 \end{align}
  \end{enumerate}
  Then $c_{ji}c_{ij} = \id$ iff $q_{ji}q_{ij} = 1$  and  $a_{ij}=0$.
  The \emph{interaction} between the block $j$ and the point $i$ is $\inc_{ij} = q_{ji}q_{ij}$. If
$q_{ij}q_{ji}= 1$,  then we say that the interaction is  weak.
  Also the \emph{ghost} between $j$ and $i$ as 
\begin{align*}
  \ghost_{ij} =  \left(-\frac{3}{2} \epsilon_j -\frac{1}{2}\right) a_{ij}.
\end{align*}
  If $\ghost_{ij} \in \N$, then we say that the ghost is \emph{discrete}.

\bigbreak
We next impose the form of the braidings between two different blocks.

  \begin{enumerate}[leftmargin=*,label=\rm{(\Alph*)}]\addtocounter{enumi}{2}
\item\label{item:braiding-2blocks} For every  $j, k \in J$, $j\neq k$, there exist
$q_{jk}, q_{kj} \in \kut$ and $a_{jk}, a_{kj} \in \ku$ such that the braiding between $V_j$ and $V_k$ with respect to the basis $B_j$ and $B_k$ as above
is given  by
\begin{align*} 
&\text{the braiding of $V_j\oplus \ku x_k$ is given by \eqref{eq:braiding-point-block} and \eqref{eq:braiding-block-point};}
\\
&\text{same for the  braiding of $V_k\oplus \ku x_j$;}
\\
&c \left(x_{\fjdos} \otimes x_{\fkdos} \right) =   q_{jk} \left(x_{\fkdos} + a_{jk} x_k \right) \otimes x_{\fjdos};
\\
&c \left(x_{\fkdos} \otimes x_{\fjdos} \right) =   q_{kj} \left(x_{\fjdos} + a_{kj} x_j \right) \otimes x_{\fkdos}.
\end{align*}
  \end{enumerate}

\bigbreak
Set   $r \sim s$ when $c_{rs}c_{sr} \neq \id_{V_s   \otimes V_r}$, $r \neq s \in I$. Let
$\approx$ be the equivalence relation on $I$ generated by $\sim$. The last assumption is: 
  
  \begin{enumerate}[leftmargin=*,label=\rm{(\Alph*)}]\addtocounter{enumi}{3}
\item\label{item:connected}
 $V$ is \emph{connected}, i.e  $r \approx s$ for all $r, s \in I$.
  \end{enumerate}
  
  \subsection{Infinite flourished graphs}\label{sec:flourished-graph}
A flourished graph is a graph $\D$ with an infinite set $\I$ of vertices and  the following decorations:
  \begin{itemize} [leftmargin=*]\renewcommand{\labelitemi}{$\circ$}
 \item The vertices have three kind of decorations $+$, $-$ and $q \in \ku^{\times}$; they are depicted respectively as $\boxplus$, $\boxminus$ and $\overset{q} {\circ}$. The set of all vertices of the first kind is denoted by $\J_+$, and those of the second kind by $\J_-$.
The vertices in $\J := \J_+ \cup \J_-$ are called blocks, the remaining are called points.

 \smallbreak
 \item If $i \neq h$ are points, and there is an edge between them, then it is decorated by some $\widetilde{q}_{ih} \in \ku^{\times} - 1$: $\xymatrix{\overset{q_i}  {\circ} \ar  @{-}[r]^{\widetilde{q}_{ih}}  & \overset{q_h}  {\circ} }$.

 \smallbreak
 \item  If $j$ is a block and $i$ is a point, then an edge between $j$ and $i$ is decorated by $\ghost_{ij}$ 
 for some $\ghost_{ij}\in \ku^{\times}$; or not decorated at all.
  \end{itemize}

The full (decorated) subgraph with vertices $\I - \J$ is denoted $\D_{\diag}$; it is a generalized Dynkin diagram \cite{H-classif} whose set of vertices
is possibly infinite.

The set of connected components of $\D_{\diag}$ is denoted by $\X$; we also set
\begin{align*}
\Xf &= \{X \in \X: \vert X \vert < \infty \},&
\Xif &= \X - \Xf. 
\end{align*}

Let $V$ be as in \S \ref{subsubsec:intro-class}. We attach a flourished graph $\D$ to $V$ by the following rules. 
The set of vertices of $\D$ is the infinite set $I$. The decoration obeys the following rules:

  \begin{itemize} [leftmargin=*]\renewcommand{\labelitemi}{$\bullet$}
 \item  If $j\in J_{+}$, respectively $j\in J_{-}$, then the  corresponding vertex is decorated as $\boxplus$, respectively $\boxminus$.
  Thus $\J_{\pm} = J_{\pm}$,  $\J = J$.
 
 \smallbreak
 \item  If $i \in I - J$,  then   the corresponding   vertex is decorated as $\overset{q_{ii}}  {\circ}$.
 
 \smallbreak
 \item There is an edge between $r$ and $s \in I$ iff $r \sim s$.
 
 \smallbreak
 \item  If $j \in J$, $i \in I - J$, $q_{ij}q_{ji} = 1$  and $a_{ij}\neq 0$, 
 then the  edge between $i$ and $j$ is labelled by $\ghost_{ij}  = \begin{cases} -2a_{ij}, & j \in J_+, \\
 a_{ij}, &j \in J_-.
 \end{cases}$.
 
 \smallbreak
 \item If $ i, h \in I - J$, $i\neq h$ and $q_{ih}q_{hi} \neq 1$, then the corresponding edge is decorated by $\widetilde{q}_{ih} = q_{ih}q_{hi}$.
  \end{itemize}
  
  \subsection{Infinite admissible graphs}\label{sec:admissible-graph}
The infinite flourished graphs arising from Nichols algebras in the class above with finite $\GK$ are described in the following definition.
  
  \begin{definition}\label{def:flourished-graph-admissible}
 An infinite flourished graph is \emph{admissible} when  the following conditions hold.
 
 \begin{enumerate}[leftmargin=*,label=\rm{(\alph*)}]
\item\label{item:finite blocks} The set $\J$ is finite and non-empty.
  
  \item\label{item:noedges} There are no edges between blocks.
  
  \smallbreak
  \item\label{item:local}
  The only possible connections between a block and a connected component $X \in \Xf$
  are described in Tables  \ref{tab:conn-comp-finitenumber} and \ref{tab:conn-comp} 
  (the point connected with the block is black for emphasis).
  Here $\ghost \in \N$, $\omega \in \G_3'$. 
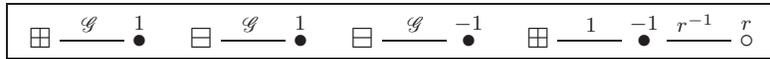
  \begin{table}[ht]
  \caption{\small{Connecting finite components and blocks;}  $r \notin \G_{\infty}$.} \label{tab:conn-comp-finitenumber}
  \begin{center}
   \begin{tabular}{|c c c c|}
   
   \hline
   $\xymatrix{\boxplus \ar  @{-}[r]^{\ghost}  & \overset{1}{\bullet}}$ &
    $\xymatrix{\boxminus \ar  @{-}[r]^{\ghost}  & \overset{1}{\bullet}}$  
   &  
    $\xymatrix{\boxminus \ar  @{-}[r]^{\ghost}  &\overset{-1}{\bullet} }$
    &  $\xymatrix{\boxplus \ar  @{-}[r]^{1} &\overset{-1}{\bullet} \ar  @{-}[r]^{r^{-1}}  & \overset{r}{\circ}}$
   \\ \hline
   \end{tabular}
  \end{center}
  
  \end{table}  
 
\begin{table}[ht]
  \caption{\small{Connecting finite components and blocks; }$r \in \G_{\infty} - \G_2$.} \label{tab:conn-comp}
  \begin{center}
 \begin{tabular}{|c  c|}
\hline
 $\xymatrix{\boxplus \ar  @{-}[r]^{\ghost}  & \overset{-1}{\bullet}}$ \quad $\xymatrix{\boxplus \ar  @{-}[r]^{1}  & \overset{\omega}{\bullet}}$
& $\xymatrix{\boxplus \ar  @{-}[r]^{1} &\overset{-1}{\bullet} \ar  @{-}[r]^{\omega^2 }  & \overset{\omega}{\circ}}$
\\ \hline  
    $\xymatrix{\overset{-1}{\circ} \ar  @{-}[r]^{\omega^2 }  & \overset{\omega}{\bullet} \ar  @{-}[r]^{1} & \boxplus }$
&  $\xymatrix{\boxplus \ar  @{-}[r]^{1}  &\overset{-1}{\bullet} \ar  @{-}[r]^{\omega}  & \overset{-1}{\circ}}$
  \\ \hline
$\xymatrix{\boxplus \ar  @{-}[r]^{1} &\overset{-1}{\bullet} \ar  @{-}[r]^{r^{-1}}  & \overset{r}{\circ}}$   
\;  $\xymatrix{\boxplus \ar  @{-}[r]^{2}  &\overset{-1}{\bullet} \ar  @{-}[r]^{-1}  & \overset{-1}{\circ}}$ 
&
$\xymatrix{\boxplus \ar  @{-}[r]^{1} &\overset{-1}{\bullet} \ar  @{-}[r]^{\omega}  & \overset{\omega^2}{\circ} \ar  @{-}[r]^{\omega}  & \overset{\omega^2}{\circ} }$
\\ \hline
$\xymatrix{\boxplus \ar  @{-}[r]^{1}  &\overset{-1}{\bullet} \ar  @{-}[r]^{-1}  & \overset{-1}{\circ}} \dots \xymatrix{
  \overset{-1}{\circ} \ar  @{-}[r]^{-1}  & \overset{-1}{\circ} }$
&  $\xymatrix{\boxplus \ar  @{-}[r]^{1} & \overset{-1}{\bullet} \ar  @{-}[r]^{\omega}  & \overset{\omega^2}{\circ} \ar  @{-}[r]^{\omega^2}  & \overset{\omega}{\circ} }$   
  \\ \hline
 \end{tabular}
  \end{center}
  
\end{table}  

\smallbreak
\item\label{item:local-finitenumber}
There are only a finite number of connections between blocks and connected components $X \in \Xf$
as in Table  \ref{tab:conn-comp-finitenumber}.

\smallbreak
\item\label{item:point} Let $X \in \Xf$. Then there is a unique $i \in X$ connected to a block.

\smallbreak
\item\label{item:concom} If $X \in \Xf$ has $\vert X \vert > 1$, then it is connected to a unique block.
  
\smallbreak
  \item\label{item:concom-1} If $X =\{i\} \in \Xf$ and $q_{ii} \in \G'_3$, then it is connected to a unique block.
  
\smallbreak
\item\label{item:graph-connected}
$\D$ is connected.  

\smallbreak
\item\label{item:local-infinite-component}

Given a connected component $X \in \Xif$, there is a unique block $V_j$ connected to $V_X$ and the corresponding flourished diagram is
\begin{align}\label{eq:f-diagram-Xif}
\xymatrix@C=12pt{\boxplus \ar  @{-}[rr]^{1}& & \overset{-1}{\bullet} \ar  @{-}[rr]^{-1}& &  \overset{-1}{\circ} \ar@{.}[rr]
  & &
  \overset{-1}{\circ} \ar  @{-}[rr]^{-1} & & \overset{-1}{\circ}  \ar@{.}[rr]
  & &  }
\end{align}

 \end{enumerate}
  \end{definition}

\begin{remark}\label{rem:cyclop}
This Definition extends \cite[Definition 1.9]{triang} to graphs with infinite sets of vertices.
Besides this, the main difference is that only weak interactions between blocks and points are allowed.
Indeed, the only possible admissible graphs in \cite[Definition 1.9]{triang} having mild interaction 
are $\cyc_1$ and $\cyc_2$, the former included in the latter, but neither contained in another admissible graph.

Another difference is that \cite[Definition 1.9]{triang} does not require connectedness but we deal with this in Corollary \ref{cor:conn-components}.
\end{remark}

\begin{remark}\label{rem:tables-3-4}
Let  $V$ be as in \S \ref{subsubsec:intro-class}; let $j \in J$, i.~e. $V_j$ is a block, and let $X \in \X$; set $V_X = \oplus_{i \in X} V_i$.
Then $\toba (V_j \oplus V_X) \simeq K \# \toba (V_j)$ for a suitable Nichols algebra $K$, see \cite[\S 4.1.4]{triang}, and
\begin{align*}
\GK \toba (V_j \oplus V_X) = \GK K + \GK \toba (V_j) = \GK K + 2.
\end{align*}
Let $\T_{\ref{tab:conn-comp-finitenumber}}$, respectively $\T_{\ref{tab:conn-comp}}$, be the set of flourished diagrams
in Table  \ref{tab:conn-comp-finitenumber}, resp.  \ref{tab:conn-comp}. 

\begin{enumerate}[leftmargin=*,label=\rm{(\alph*)}]
  \item\label{item:table3} If the diagram of $V_j \oplus V_X$ belongs to $\T_{\ref{tab:conn-comp-finitenumber}}$, then
$\GK \toba (V_j \oplus V_X) \geq 3$. 

  \item\label{item:table4} If the diagram of $V_j \oplus V_X$ belongs to $\T_{\ref{tab:conn-comp}}$, then
  $\GK \toba (V_j \oplus V_X) = 2$.
\end{enumerate}

See \cite[Tables 2 and 3]{triang}, and references therein.
\end{remark}

\begin{theorem}\label{theorem:main} 
Let $V$ be a braided vector space as in \S \ref{subsubsec:intro-class} and let $\D$ be its infinite flourished graph. 
The following are equivalent:

\begin{enumerate}[leftmargin=*,label=\rm{(\Roman*)}]
\item\label{item:finiteGK} $\GK \toba(V) < \infty$,

\item\label{item:admissible} $\D$ is admissible.
 \end{enumerate}
\end{theorem}  

\pf \ref{item:finiteGK} $\implies$ \ref{item:admissible}: First,
\ref{item:noedges} follows from \cite[Theorem 6.1]{triang}. Now $J \neq \emptyset$ in
\ref{item:finite blocks} is part of the assumption \ref{item:decomposition}. Let $j_1, \dots, j_t$ be different blocks.
Then $\GK \toba(V_{j_1} \oplus \dots \oplus V_{j_t}) = 2t$ by the proof of \cite[Theorem 7.1]{triang}.
Hence $J$ is finite.

Let $j \in J$ be a block and $X \in \Xf$ connected to $j$. Then the interaction between them is weak as explained in Remark \ref{rem:cyclop}. 
By \cite[Theorem 1.10]{triang}, \ref{item:local}, \ref{item:point}, \ref{item:concom} and \ref{item:concom-1} follow.

Let $\cV_1 = \oplus_{j\in J} V_j$ and let $X_1, \dots, X_m \in \Xf$ be such that the connection between $X_l$ and a block  is as in Table \ref{tab:conn-comp-finitenumber},
for every $l \in \I_m$. Let $\cV_2 = \oplus_{i\in X_1 \cup \dots \cup X_m} V_i$.
Then
\begin{align*}
\GK \toba (V) \geq \GK \toba (\cV_1  \oplus \cV_2) \geq 2\vert J\vert + m,
\end{align*}
by the formula at the end of the proof of  \cite[Theorem 7.1]{triang}, together with Remark \ref{rem:tables-3-4}.
This shows \ref{item:local-finitenumber}. 

Also,   \ref{item:graph-connected} is the assumption \ref{item:connected}.
Finally, if $X \in \Xif$, then  it is connected to a block $j$ by  \ref{item:connected}. 
Then \ref{item:local} and \ref{item:concom} 
say that $X$ and $j$ should have the form in \ref{item:local-infinite-component}.

\ref{item:admissible} $\implies$ \ref{item:finiteGK}: By \ref{item:local}, we have a splitting $\Xf = \X_{\ref{tab:conn-comp-finitenumber}} \coprod \X_{\ref{tab:conn-comp}}$ where 
\begin{align*}
\X_{\ref{tab:conn-comp-finitenumber}} = &\{X \in \Xf: \exists j\in J \text{ such that } V_j \oplus V_X \text{ has diagram in } \T_{\ref{tab:conn-comp-finitenumber}}\},
\\
\X_{\ref{tab:conn-comp}} = &\{X \in \Xf: \exists j\in J \text{ such that } V_j \oplus V_X \text{ has diagram in } \T_{\ref{tab:conn-comp}}\}.
\end{align*}
By \ref{item:finite blocks} and \ref{item:local-finitenumber}, the braided vector subspace
\begin{align*}
\cV_0 &= \left(\oplus_{j\in J} V_j\right) \oplus \left(\oplus_{X\in \X_{\ref{tab:conn-comp-finitenumber}}} V_X\right)
\end{align*}
has finite dimension. By \cite[Theorem 7.1]{triang}, cf. Remark \ref{rem:cyclop},  
\begin{align*}
d &:= \GK \toba(\cV_0)  < \infty.
\end{align*}
Given $Y \in \Xif$ and $n \in \N$, we denote by $Y[n]$ the connected subdiagram of $Y$ with $n$ vertices starting at the black point.
Let us now consider finite subsets $F \subset \X_{\ref{tab:conn-comp}}$ and $G \subset \Xif$, together with a function $\mathbf{n}: G\to \N$, 
$Y \mapsto n_Y$. We set
\begin{align*}
\cV_{F, G, \mathbf{n}} &= \cV_0 \oplus \left(\oplus_{X\in F} V_X\right) \oplus \left(\oplus_{Y\in G} V_{Y[n_Y]}\right).
\end{align*}
By the proof of \cite[Theorem 7.1]{triang}, cf. Remark \ref{rem:cyclop},  
\begin{align*}
\GK \toba(\cV_{F, G, \mathbf{n}}) &=  d.
\end{align*}
Since $V$ is the filtered union of all the $\cV_{F, G, \mathbf{n}}$'s, we conclude by Lemma \ref{prop:infinite-dynkin-rels}
that $\GK \toba(V) =  d$.   \epf

Let now $V$ be a braided vector space as in \ref{subsubsec:intro-class} except that we do not assume \ref{item:connected}, 
i.~e. connectedness. Let $\Zh$ be the set of connected components of $V$ (do not confuse with the set $\X$ of connected components of $V_{\diag}$). Given $K \subset I$, we set as above $V_{K} = \oplus_{i \in K} V_i$. Let
\begin{align*}
\Zh_{>0} &= \{\zh \in \Zh: \GK \toba(V_{\zh}) > 0 \}.
\end{align*}

\begin{lemma}\label{lemma:GK-dim-sum}
Let $I_1$ be a proper non-empty subset of $I$ and $I_2 = I - I_1$. 
If $c_{hi}c_{ih} = \id_{V_i \otimes V_h}$ for all $i \in I_1$ and $h \in I_2$, then
\begin{align*}
\GK \toba(V) = \GK \toba(V_{I_1}) + \GK \toba(V_{I_2}).  
\end{align*}
\end{lemma}

\pf We may assume that $\GK \toba(V_{I_1}) < \infty$ and  $\GK \toba(V_{I_2}) < \infty$.
Let $F$ be a finite subset of $I$ and $F_a = F \cap I_a$, $a = 1,2$, thus $F = F_1 \cup F_2$. 
Then $\GK \toba(V_F) = \GK \toba(V_{F_1}) + \GK \toba(V_{F_2})$ since $\toba (V_F) \simeq \toba (V_{F_1}) \underline{\otimes} \toba (V_{F_2})$
and both have convex PBW-basis, hence GK-deterministic subspaces, see Remark \ref{rem:GK-det} and \cite[Lemma 3.1]{triang}. Hence Lemma \ref{prop:infinite-dynkin-rels} applies.  
\epf

\begin{coro}\label{cor:conn-components}
The following are equivalent:
\begin{enumerate}[leftmargin=*,label=\rm{(\Roman*)}]
\item\label{item:finiteGK-cor} $\GK \toba(V) < \infty$.
  
\item\label{item:admissible-cor} $\Zh_{>0}$ is finite; and
for each $\zh \in \Zh$, $\GK \toba(V_{\zh}) < \infty$, either $V_{\zh}$ is of diagonal type or else it has an admissible flourished diagram.

\end{enumerate}
\end{coro}

\pf  \ref{item:finiteGK} $\implies$ \ref{item:admissible}: If $\zh_1, \dots, \zh_d$ are different components in 
$\Zh_{>0}$, then $\GK \toba(V) \ge d$ by Lemma \ref{lemma:GK-dim-sum}. The second  statement is evident 
and the third  follows from Theorem \ref{theorem:main}. 

\ref{item:admissible} $\implies$ \ref{item:finiteGK}: By Lemma \ref{lemma:GK-dim-sum},  
$\GK \toba(\oplus_{\zh \in \Zh_{>0}} V_{\zh}) < \infty$; call it $d$. Then 
$\GK \toba\left(\oplus_{\zh \in F} V_{\zh}\right) < \infty$
for any finite subset $F$ of $\Zh$ that contains $\Zh_{>0}$ by the same result.
By Lemma \ref{prop:infinite-dynkin-rels}, the claim follows.
\epf

\subsection{Examples}\label{subsec:examples}
We illustrate the previous result describing some examples of Nichols algebras of infinite rank and finite $\GK$.

\begin{example} Let $\Iw = \N \cup \{\frac{3}{2} \}$.
Let $\lstr(A_{\infty})$ be the braided vector space defined by a matrix $(q_{ij})_{i,j\in \N}$ in such a way that it has a flourished diagram 
\begin{align*}
\xymatrix{\underset{1}{\boxplus} \ar  @{-}[rr]^{1}& & \overset{-1}{\underset{2}{\bullet}} \ar  @{-}[rr]^{-1}& &  \overset{-1}{\underset{3}{\circ}} \ar@{.}[rr]
& &
\overset{-1}{\underset{j}{\circ}} \ar  @{-}[rr]^{-1} & & \overset{-1}{\underset{j+1}{\circ}}  \ar@{.}[r]
&  }
\end{align*}
By Corollary \ref{cor:infinite-dynkin-rels} and \cite[Proposition 5.31]{triang},
the algebra $\toba(\lstr(A_{\infty}))$ has $\GK = 2$.
Also it is presented by generators $x_i$, $i\in \Iw_{\theta}$ with relations
as in \cite[Proposition 5.31]{triang},  replacing $\theta$ by $\infty$.
A PBW basis is obtained by union of PBW-basis of the algebras $\toba(\lstr(A_{\theta}))$, $\theta \in \N$.
 \end{example}

\begin{example} Let $(n_k)_{k \in \N_{\ge 2}}$ be a family of natural numbers and
  $\Iw = \bigcup_{k \in \N_{\ge 2}}(\{k\} \times \I_{n_k} ) \cup \{1, \frac{3}{2} \}$.
  Let $V$ be the braided vector space with flourished diagram 
  \begin{align*}
  \xymatrix@R=0.5pt{& & \overset{-1}{\underset{(2, 1)}{\bullet}} \ar  @{-}[rr]^{-1}& &  \overset{-1}{\underset{(2,2)}{\circ}} \ar@{.}[rr]
    & &
    \overset{-1}{\underset{(2, n_2)}{\circ}}  
    &  
\\
  & & \overset{-1}{\underset{(3,1)}{\bullet}} \ar  @{-}[rr]^{-1}& &  \overset{-1}{\underset{(3,2)}{\circ}} \ar@{.}[rr]
& &
\overset{-1}{\underset{(3, n_3)}{\circ}}  
&  
\\
\underset{1}{\boxplus} \ar  @{-}[uurr]^{1} \ar  @{-}[urr]^{1}\ar  @{-}[drr]^{1} \ar @{.}[ddrr]    \ar  @{.}[rrrrrr] &&  && &&
\\
& & \overset{-1}{\underset{(k, 1)}{\bullet}} \ar  @{-}[rr]^{-1}& &  \overset{-1}{\underset{(k,2)}{\circ}} \ar@{.}[rr]
& &
\overset{-1}{\underset{(k, n_k)}{\circ}}    
\\
&& \ar  @{.}[rrrr]  && && }
  \end{align*}
By Corollary \ref{cor:infinite-dynkin-rels} and \cite[Proposition 5.31]{triang}
the algebra $\toba(V)$  is presented by generators $x_i$, $i\in \Iw_{\theta}$, 
with the relations of the  various subalgebras $\toba(\lstr(A_{n_k -1}))$ together with 
$q$-commuting relations between the points in different $A_{n_k -1}$'s (but with various $q$'s). 
It has $\GK = 2$ and a PBW-basis is constructed along the lines of the proof of 
 \cite[Theorem 7.1]{triang}.
 
Variation: replace some (or all) the $n_k$'s by $\infty$.
\end{example}

\begin{example} 
Let $\Iw = \N \cup \{\frac{3}{2}, \frac{5}{2} \}$. Let $(\ghost_{i1})_{i\in \N_{\ge 3}}$, $(\ghost_{i2})_{i\in \N_{\ge 3}}$
  be two families of natural numbers and $\bq = (q_{ij})_{i, j \in \N}$ giving rise to the flourished diagram 
  
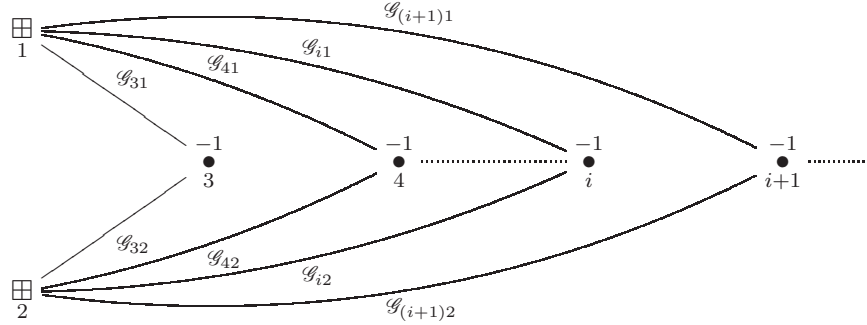
\begin{align*}
\xymatrix{\underset{1}{\boxplus} \ar  @{-}[drr]^{\ghost_{31}} \ar  @/^0.5pc/@{-}[drrrr]^{\ghost_{41}}  \ar  @/^1pc/@{-}[drrrrrr]^{\ghost_{i1}} \ar  @/^2pc/@{-}[drrrrrrrr]^{\ghost_{(i+1)1}}  & &   & & & &   & & &
\\ & & \overset{-1}{\underset{3}{\bullet}}  & &  \overset{-1}{\underset{4}{\bullet}} \ar@{.}[rr]
    & &
    \overset{-1}{\underset{i}{\bullet}}   & & \overset{-1}{\underset{i+1}{\bullet}}  \ar@{.}[r]
    &  
\\
\underset{2}{\boxplus} \ar  @{-}[urr]_{\ghost_{32}} \ar  @/_0.5pc/@{-}[urrrr]_{\ghost_{42}}  \ar  @/_1pc/@{-}[urrrrrr]_{\ghost_{i2}} \ar  @/_2pc/@{-}[urrrrrrrr]_{\ghost_{(i+1)2}} & &   & & & &   & & &  
}
  \end{align*}
Let $V$ be the braided vector space with this diagram; notice that the subdiagram spanned by $\{1,2,i\}$ 
corresponds to a Poseidon braided subspace $\pos_i$, as in \cite[\S 7]{triang}, for every
$i\in \N_{\ge 3}$.
By Corollary \ref{cor:infinite-dynkin-rels},
the algebra $\toba(V)$ is presented by generators $x_i$, $i\in \Iw$, with the defining relations
of the various $\toba(\pos_i)$, cf. \cite[Proposition 7.7]{triang}, together with the $q_{ih}$ -commuting relations for 
$i\neq h\in \N_{\ge 3}$.
It has $\GK = 4$ and a PBW-basis by collecting together those of the  various $\toba(\pos_i)$, 
cf. the proof of  \cite[Theorem 7.1]{triang}.
\end{example}

\subsection{Hopf algebras with finite $\GK$}\label{subsec:examples-hopf}
Let $V$ be a braided vector space as in \S \ref{subsubsec:intro-class}; assume that its flourished diagram
is admissible. 

A \emph{principal realization} of $V$ over an abelian group $\Gamma$ consists of

\begin{enumerate}[leftmargin=*,label=\rm{(\roman*)}]
\item a family $(g_i)_{i\in I}$ of elements of $\Gamma$,
\item a family $(\chi_i)_{i\in I}$ of characters of $\Gamma$,
\item a family $(\eta_j)_{j\in J}$ of derivations of $\Gamma$,    
\end{enumerate}
such that
\begin{align}
\label{eq:realization-qij}
\chi_h(g_i) &= q_{ih}, & i,h &\in I,
\\
\eta_j(g_i) &= a_{ij}, & i &\in I, \, j\in J.
\label{eq:realization-derivation}
\end{align}

Given a principal realization the braided vector space $V$ is realized in $\ydG$, hence we get a Hopf algebra
by bosonization $\toba(V) \# \ku \Gamma$. Notice that the realization depends not only on the Dynkin diagram but actually on all the $q_{ij}$'s.
For convenient choices of the last, one can find an abelian group $\Gamma$ which is finitely generated modulo its torsion.
Then $\GK \toba(V) \# \ku \Gamma$ would be finite. We leave to the reader the exercise of working out these ideas.

\end{document}